\documentclass[10pt]{article}
\usepackage{fancyhdr,a4wide}
\usepackage[mathscr]{eucal}
\usepackage{amsmath}
\usepackage{amssymb}
\usepackage{amsthm}

\usepackage{caption}
\usepackage{enumerate}
\usepackage{fullpage}
\usepackage{graphicx}
\usepackage{multicol}
\usepackage{tikz, subcaption, mathtools, soul, mathdots, yhmath}
\usetikzlibrary{calc,decorations.markings,arrows}

\usepackage{hyperref}

\definecolor{background}{RGB}{230,230,230} 
\definecolor{arc1}{RGB}{0,128,128}   
\definecolor{arc2}{RGB}{210,105,30}	
\definecolor{arc3}{RGB}{80,0,80}		
\definecolor{arc4}{RGB}{34,139,34}		

\input xy
\xyoption{all}

\usepackage[all]{xy}

\title{Frieze patterns of integers} 
\author{Karin Baur} 

\usepackage{pst-all}  

\definecolor{light-grey}{gray}{0.6}  

\begin{document}

\newtheorem{lm}{Lemma}[section]
\newtheorem{prop}[lm]{Proposition}
\newtheorem{satz}[lm]{Satz}

\newtheorem*{corollary}{Corollary}
\newtheorem{cor}[lm]{Korollar}
\newtheorem{theorem}[lm]{Theorem}
\newtheorem*{thm}{Theorem}

\theoremstyle{definition}
\newtheorem{defn}[lm]{Definition}
\newtheorem{qu}{Question}
\newtheorem{ex}[lm]{Example}
\newtheorem{exas}[lm]{Examples}
\newtheorem{exc}[lm]{Exercise}
\newtheorem*{facts}{Facts}
\newtheorem{rem}[lm]{Remark}

\theoremstyle{remark}


\newcommand{\perm}{\operatorname{Perm}\nolimits}
\newcommand{\NN}{\operatorname{\mathbb{N}}\nolimits}
\newcommand{\ZZ}{\operatorname{\mathbb{Z}}\nolimits}
\newcommand{\QQ}{\operatorname{\mathbb{Q}}\nolimits}
\newcommand{\RR}{\operatorname{\mathbb{R}}\nolimits}
\newcommand{\CC}{\operatorname{\mathbb{C}}\nolimits}
\newcommand{\PP}{\operatorname{\mathbb{P}}\nolimits}
\newcommand{\cA}{\operatorname{\mathcal{A}}\nolimits}
\newcommand{\LL}{\operatorname{\Lambda}\nolimits}

\newcommand{\MM}{\operatorname{\mathcal M}\nolimits}
\newcommand{\FF}{\operatorname{\mathcal F}\nolimits}

\newcommand{\Fk}{\mathcal{F}_{k,n}}
\newcommand{\Mk}{M_{k,n}}

\newcommand{\ad}{\operatorname{ad}\nolimits}
\newcommand{\im}{\operatorname{im}\nolimits}
\newcommand{\Char}{\operatorname{char}\nolimits}
\newcommand{\Aut}{\operatorname{Aut}\nolimits}

\newcommand{\id}{\operatorname{id}\nolimits}
\newcommand{\ord}{\operatorname{ord}\nolimits}
\newcommand{\ggT}{\operatorname{ggT}\nolimits}
\newcommand{\lcm}{\operatorname{lcm}\nolimits}

\newcommand{\Gr}{\operatorname{Gr}\nolimits}

\maketitle

\begin{abstract}
The famous theorem of Conway and Coxeter on frieze patterns gave a geometric 
interpretation to integral friezes via triangulations of polygons. 
In this article, we review 
this result and show some of the development it has led to. 
The last decade has seen a lot of activities on friezes. 
One reason behind this is the connection to cluster combinatorics. 
\end{abstract}

%
%
\section{Introduction}\label{sec:intro}
%

The first time I encountered John Conway's work on frieze patterns was around 2006 when I participated 
in a a reading group on 
cluster algebras and related topics in the Pure Mathematics group at the University of Leicester. 
The two papers~\cite{Coco1, Coco2} were on our reading list at the time. These two articles are very remarkable as they 
invite the reader to a discovery of the properties of frieze patterns. The first of them presents 
a list of 39 problems to study whereas the second paper gives hints and solutions for these problems. 
Problems (28) and (29) are known as the theorem of Conway and Coxeter (Theorem~\ref{thm:CC} below), relating 
integer friezes with triangulations of polygons: 

\noindent
(28) 
Is there just one frieze pattern of integers for each triangulated polygon? 

\noindent 
(29)
Is there just one triangulated polygon for each frieze pattern of integers? 

These two problems have lead to a huge amount of activities around friezes and in this article we 
give an overview 
of the development initiated by Conway's contribution to friezes almost 50 years ago. It will 
be by no means exhaustive, there is not enough space here to give credit to all the 
development in this area. 
There are currently 88 articles on mathscinet with the word ``frieze'' in the title. Some of them 
appeared in the 20 years since their inception. However, the 
discovery of the link between friezes and cluster algebras by Fomin and Zelevinsky, see~\cite{FZ03} 
and between friezes and cluster categories by Caldero and Chapoton, \cite{CCh06} 
led to renewed interest in this topic, 
with over 60 articles being written around friezes in the last 10 years. 
The original definition of Coxeter, \cite{C71} 
has been generalised in many directions, 
resulting in work on infinite friezes, SL$_k$-friezes, SL$_k$-tilings, tropical friezes, 2-friezes, 
$q$-deformed friezes (as the article \cite{MGO20} in this volume), etc. 
It became apparent that friezes are of interest to various areas such as combinatorics, 
cluster algebras, geometry, 
integrable systems and representation theory. 
A good survey of this is the article~\cite{Crossroads} by Sophie Morier-Genoud where many of these 
directions and links to the literature can be found. 
Frieze patterns were the topic of a mini-workshop in 
Oberwolfach in 2015, \cite{mfofriezes}; an MFO snapshot authored by Thorsten Holm appeared 
subsequently, \cite{Holm2015}. 
The combinatorial flavour of friezes have them made very popular for outreach activities. 
There are even two beautiful videos on friezes by Numberphile, presented 
by Sergei Tabachnikov, \cite{friezevideo1}, \cite{friezevideo2}. 

In this article, we concentrate on the combinatorical geometric interpretation of friezes inspired 
by problems (28) and (29) above. 
We first recall the notion of a frieze pattern (Section~\ref{sec:friezes}) and some of the 
properties of closed integer friezes. Then we concentrate on infinite friezes in Section~\ref{sec:infinite}, 
giving a geometric interpretation via matchings and explaining the growth of the entries in 
those friezes. 

\subsection*{Acknowledgement}

The author thanks Sophie Morier-Genoud for helpful comments and suggestions. 
The author was supported by Royal Society Wolfson Fellowship 180004 
and by FWF projects P 30549.  She acknowledges support by FIM, ETHZ. 
She is currently on leave from the University of Graz. 

%
\section{Frieze patterns}\label{sec:friezes}
%

Frieze patterns have first been studied by H.S.M. Coxeter in~\cite{C71}. 
We recall their definition and then 
discuss finite and infinite integral frieze patterns.

\begin{defn}\label{def:frieze}
A {\bf frieze pattern} is 
formed by a finite or infinite number of rows, shifted with respect to each other, starting with a row of 0s, 
followed by a row of 1s and then rows of elements $a_{ij}$, $i+2\le j$, of an integral domain 
(mostly the integral domain will be the integers here - but there will be an instance where the 
entries are cluster variables of a cluster algebra) : 

\begin{center}
\begin{tikzpicture}[scale=.8, inner sep=4pt]
\node (00) at (-6,6) {$\ldots$};
\node (01) at (-4,6) {0};
\node (02) at (-2,6) {0};
\node (ai+) at (0,6) {0};
\node (03) at (2,6) {0};
\node (04) at (4,6) {0};
\node (05) at (6,6) {$\ldots$};
\node (10) at (-7,5) {$\ldots$};
\node (11) at (-5,5) {1};
\node (12) at (-3,5) {1};
\node (12) at (-1,5) {1};
\node (12) at (1,5) {1};
\node (14) at (3,5) {1};
\node (15) at (5,5) {1};
\node (16) at (7,5) {$\ldots$};

\node (a-2) at (-6,4) {$\ldots$};
\node (a-1) at (-4,4) {\(a_{-1,1}\)};
\node (a0) at (-2,4) {\(a_{0,2}\)};
\node (a1) at (0,4) {\(a_{1,3}\)};
\node (a2) at (2,4) {\(a_{2,4}\)};
\node (a3) at (4,4) {\(a_{3,5}\)};
\node (a4) at (6,4) {$\ldots$};

\node (10) at (-7,3) {$\ldots$};
\node (11) at (-5,3) {\(a_{-2,1}\)};
\node (12) at (-3,3) {\(a_{-1,2}\)};
\node (12) at (-1,3) {\(a_{0,3}\)};
\node (12) at (1,3) {\(a_{1,4}\)};
\node (14) at (3,3) {\(a_{2,5}\)};
\node (15) at (5,3) {\(a_{3,6}\)};
\node (16) at (7,3) {$\ldots$};

\node (ai+) at (-6,2) {$\ldots$};
\node (ai) at (-4,2) {\(a_{-2,2}\)};
\node (ai+) at (-2,2) {\(a_{-1,3}\)};
\node (ai+) at (0,2) {\(a_{0,4}\)};
\node (aj-) at (2,2) {\(a_{1,5}\)};
\node (aj) at (4,2) {\(a_{2,6}\)};
\node (ai+) at (6,2) {$\ldots$};

\node (ai) at (-5,1) {$\vdots$};
\node (ai+) at (-3,1) {$\vdots$};
\node (ai+) at (-1,1) {$\vdots$};
\node (aj-) at (1,1) {$\vdots$};
\node (aj) at (3,1) {$\vdots$};
\node (aj) at (5,1) {$\vdots$};

\end{tikzpicture}
\end{center}

There are two conditions on these patterns: we ask first that 
any square formed by four neighboured entries satisfies the {\bf unimodular} or {\bf diamond rule}: 
we have $bc-ad=1$ for any four entries
\[
\xymatrix@C=.3em@R=.2em{ & a \\ b && c \\ & d} 
\]
Furthermore, we ask that such a pattern is periodic: there is an $n>0$ such that $a_{i,j}=a_{i+n,j+n}$ for any 
$i+2\le j$. 
\end{defn}

By the above definition, if the entries $a_{ij}$ are all positive integers, the frieze pattern 
is determined by its first non-trivial row. This is true more generally (for integer domains) if we 
ask that every $3\times 3$ diamond in the frieze has determinant $0$. Such frieze pattern are 
said to be {\bf tame}. 
If the frieze is 
$n$-periodic, any $n$ successive entries $a_{i,i+2},\dots, a_{i+n-1,i+n+1}$ determine the frieze under 
iterated application of the diamond rule. Such a tuple is a {\bf quiddity sequence} for the frieze. 

A frieze pattern is called {\bf integral} if all entries $a_{ij}$ are positive integers. It is called {\bf closed} if 
after a finite number of rows it stops with a second row of 1s followed by a row of 0s. Otherwise, the frieze pattern is 
{\bf infinite}. 
The {\bf order} of closed  frieze is defined to by one less than its number of 
rows (including the 0s). (The number of non-trivial rows is called the width of the frieze.) 
An example of a closed integral frieze of order 6 is in Figure~\ref{fig:order6}. 
An example of an infinite frieze is on the right in Figure~\ref{fig:infinite}. 

\begin{figure}[ht!]
\[
\xymatrix@C=1em@R=.6em{
0 && 0 && 0 && 0 && 0 && 0 && 0 && \\
& 1 && 1 && 1 && 1 && 1 && 1 &&  \\
2 && 1 && 4 && 1 && 2 && 2 && 2 &&  \\
& 1 && 3 && 3 && 1 && 3 && 3 &&  \\
1 && 2 && 2 && 2 && 1 && 4 && 1 && \\
& 1 && 1 && 1 && 1 && 1 && 1 &&  \\
0 && 0 && 0 && 0 && 0 && 0 && 0 && 
}
\]
\caption{A closed integral frieze of order 6}\label{fig:order6}
\end{figure}

From now on we will mostly concentrate on integral friezes. In this case, we have a geometric interpretation 
of frieze patterns via triangulations of polygons (if the frieze pattern is closed) or of annuli (if it is infinite) as 
we will see. 

%
\subsection{Closed frieze patterns}

First note that if $\mathcal F$ is a closed integral frieze of order $n$, then $\mathcal F$ is 
$n$-periodic, cf.~\cite[(21)]{Coco1}. 
This can be viewed as a consequence of the Conway-Coxeter theorem: 

\begin{theorem}\cite[Problems (28), (29)]{Coco1,Coco2}\label{thm:CC}
There is a bijection between frieze patterns of order $n$ and triangulations of convex $n$-gons. 
\end{theorem}

We will now explain how this works via matchings for triangulations. 
Consider a convex polygon $P_n$ with $n$ vertices, labeled 
clockwise by $\{1,2,\dots,n\}$. Take a triangulation of $P_n$, i.e. a maximal collection by pairwise 
non-crossing diagonals of $P_n$ (any such collection contains $n-3$ diagonals). 
Let $a_i$ be the number matchings of triangles with 
vertex $i$, i.e.\ the number of triangles incident with vertex $i$. 
Then $(a_1,\dots, a_n)$ is the quiddity sequence of a frieze pattern of $n$. 

Some of the first properties of frieze patterns are that every quiddity sequence 
of a closed integral frieze contains at least two entries equal to 1: The above theorem gives 
a geometric reason for it  
since every triangulation contains at least two triangles which are peripheral, i.e.~where 
two edges are boundary segments. Another fact is that if $(a_1,\dots, a_n)$ is a quiddity sequence 
with $n>2$, there can be no two entries 1 next to each other. 

\begin{rem}\label{rem:cut-glue}
There is a way to go from quiddity sequences with $n$ elements to quiddity sequences with 
$n+1$ elements and back, as described in \cite[(23)]{Coco1}: If 
$(a_1,\dots, a_n)$ is the quiddity sequence of a frieze pattern of order $n$ and $i$ some 
index, $1\le i\le n$, the 
sequence $(a_1,\dots, a_{i}+1,1,a_{i+1}+1,\dots, a_n)$ is the quiddity sequence of a 
frieze pattern of order 
$n+1$. We call this operation ``gluing''. 
The reverse operation to gluing is to start with an entry $a_i=1$ in a quiddity sequence and 
to remove it, by simultaneously decreasing its two neighbours by 1: from 
$(a_1,\dots, a_{i-1},1, a_{i+1},\dots, a_n)$ we obtain the new quiddity sequence 
$(a_1,\dots, a_{i-1}-1,a_{i+1}-1,\dots, a_n)$ with $n-1$ entries. This is called ``cutting (at an 
entry 1)''. 

In terms of triangulations of polygons, gluing corresponds to adding a vertex with a 
peripheral triangle to the triangulated polygon, the operation cutting 
to removing a peripheral triangle. 
\end{rem}

\begin{rem}[Idea of the proof of Theorem~\ref{thm:CC}]
One checks that the claim is true for $n=3$ and then uses induction. 
To go from a frieze pattern $\FF$ of order $n\ge 2$ to a triangulated $n$-gon, one modifies the 
quiddity sequence of $\FF$ by cutting at an entry 1. The result is a frieze pattern $\FF'$ of 
order $n-1$. By induction hypothesis, this corresponds to a triangulation of a polygon on 
$n-1$ vertices. Glue a triangle at the appropriate position to get the triangulation searched for. 

To go from a triangulated $n$-gon with $n\ge 4$ to a frieze, one uses the fact that every triangulated 
polygon has a peripheral triangle. Removing this gives a triangulated polygon with $n-1$ vertices. 
This corresponds to a closed frieze of order $n-1$, by induction. 
Extending its quiddity sequence by gluing yields 
the quiddity sequence for the triangulated $n$-gon we wanted. 
\end{rem}

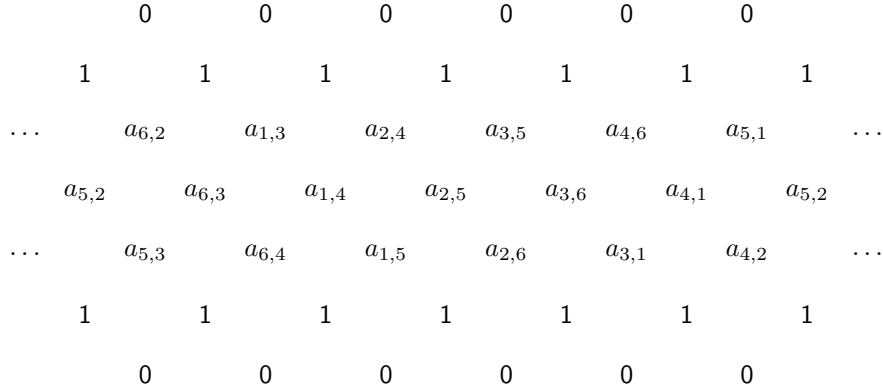
\begin{figure}[ht!]
\begin{center}
\begin{tikzpicture}[scale=.8, inner sep=4pt]
\node (01) at (-4,6) {0};
\node (02) at (-2,6) {0};
\node (ai+) at (0,6) {0};
\node (03) at (2,6) {0};
\node (04) at (4,6) {0};
\node (04) at (6,6) {0};

\node (11) at (-5,5) {1};
\node (12) at (-3,5) {1};
\node (12) at (-1,5) {1};
\node (12) at (1,5) {1};
\node (14) at (3,5) {1};
\node (15) at (5,5) {1};
\node (16) at (7,5) {1};

\node (a-2) at (-6,4) {$\ldots$};
\node (a-1) at (-4,4) {\(a_{6,2}\)};
\node (a0) at (-2,4) {\(a_{1,3}\)};
\node (a1) at (0,4) {\(a_{2,4}\)};
\node (a2) at (2,4) {\(a_{3,5}\)};
\node (a3) at (4,4) {\(a_{4,6}\)};
\node (a4) at (6,4) {\(a_{5,1}\)};
\node (16) at (8,4) {$\ldots$};

\node (11) at (-5,3) {\(a_{5,2}\)};
\node (12) at (-3,3) {\(a_{6,3}\)};
\node (12) at (-1,3) {\(a_{1,4}\)};
\node (12) at (1,3) {\(a_{2,5}\)};
\node (14) at (3,3) {\(a_{3,6}\)};
\node (15) at (5,3) {\(a_{4,1}\)};
\node (16) at (7,3) {\(a_{5,2}\)};


\node (ai+) at (-6,2) {$\ldots$};
\node (ai) at (-4,2) {\(a_{5,3}\)};
\node (ai+) at (-2,2) {\(a_{6,4}\)};
\node (ai+) at (0,2) {\(a_{1,5}\)};
\node (aj-) at (2,2) {\(a_{2,6}\)};
\node (aj) at (4,2) {\(a_{3,1}\)};
\node (aj+) at (6,2) {\(a_{4,2}\)};
\node (ai++) at (8,2) {$\ldots$};

\node (ai) at (-5,1) {1};
\node (ai+) at (-3,1) {1};
\node (ai+) at (-1,1) {1};
\node (aj-) at (1,1) {1};
\node (aj) at (3,1) {1};
\node (aj) at (5,1) {1};
\node  at (7,1) {1};

\node (ai) at (-4,0) {0};
\node (ai+) at (-2,0) {0};
\node (ai+) at (0,0) {0};
\node (aj-) at (2,0) {0};
\node (aj) at (4,0) {0};
\node  at (6,0) {0};

\end{tikzpicture}
\end{center}
\caption{The layout of a frieze of order 6.} 
\label{fig:order6-labels}
\end{figure}

By Theorem~\ref{thm:CC} we can give the entries of any closed integral frieze labels of the vertices of the polygon. 
The entries in the quiddity sequence $(a_{n,2},a_{1,3},\dots, a_{n-1,1}$) are viewed to be in positions 
$(n,2), (1,3)$, $(2,4)$, $\dots, (n-1,1)$ and going along a diagonal NW-SE, we keep the first coordinate, 
the entries in diagonals NE-SW have the second entry fixed. 
Figure~\ref{fig:order6-labels} shows this for $n=6$. 

With this labelling system we can explain how every entry in a closed integral frieze is 
given by a set of matchings for a triangulation, a result of Broline, Crowe and Isaacs. 

Consider a triangulated convex $n$-gon $P$. Let $i,j$ be two vertices of $P$ with 
$j\notin \{i-1,i,i+1\}$. 
Denote by $\MM(i,j)$ the set of matchings of triangles with the set of vertices 
$\{i+1,\dots, j-1\}$ (reducing modulo $n$). Note that 
$\{i+1,\dots, j-1\}\cup  \{j+1,\dots, i-1\}=\{1,\dots,n\}\setminus \{i,j\}$. 

\begin{theorem}\cite[Theorem 1]{BCI74}\label{thm:bci-1} 
Let $T$ be a triangulation of a convex $n$-gon $P$ and let $\FF=(a_{ij})_{ij}$ 
be the frieze associated to $T$. 
Then $|\MM(i,j)|=a_{ij}=|\MM(j,i)|$ for any two vertices $i,j$ of $P$ with $j\notin \{i-1,i,i+1\}$. 
\end{theorem}

\begin{ex}\label{ex:6gon}
Consider a fan triangulation of a hexagon with all diagonals incident with vertex $1$. 

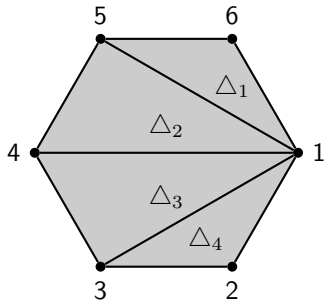
\begin{minipage}{.4\textwidth}
\begin{tikzpicture}[scale=0.35]
    \fill[gray!40] (0:5)-- (300:5) -- (240:5) -- (180:5)  -- (120:5) -- (60:5) -- cycle;

    \node[label=right:1,shape=circle,fill=black,scale=.4] (o1) at (0:5) {}; 
    \node[label=below:2,shape=circle,fill=black,scale=.4] (o2) at (300:5) {};
    \node[label=below:3,shape=circle,fill=black,scale=.4] (o3) at (240:5) {};
    \node[label=left:4,shape=circle,fill=black,scale=.4] (o4) at (180:5) {};
    \node[label=above:5,shape=circle,fill=black,scale=.4] (o5) at (120:5) {};
    \node[label=above:6,shape=circle,fill=black,scale=.4] (o6) at (60:5) {};

    
   \draw[thick, black] (o1) to (o2); 
   \draw[thick, black] (o3) to (o2);
   \draw[thick, black] (o3) to (o4);
   \draw[thick, black] (o5) to (o4);
   \draw[thick, black] (o5) to (o6);
   \draw[thick, black] (o1) to (o6);
   \draw[thick, black] (o1) to (o3);
   \draw[thick, black] (o1) to (o4);
   \draw[thick, black] (o1) to (o5);
   
   \draw (0,1) node     {$\triangle_2$};
   \draw (2.5,2.5) node     {$\triangle_1$};
   \draw (0,-1.7) node     {$\triangle_3$};
   \draw (1.5,-3.3) node     {$\triangle_4$};
\end{tikzpicture}
\end{minipage}
\begin{minipage}{.55\textwidth}
The quiddity sequence of the corresponding frieze is $(4,1,2,2,2,1)$, as in Figure~\ref{fig:order6}. 
To illustrate Theorem~\ref{thm:bci-1}, 
we determine a few sets of matchings. 
There are four triangles incident with vertex $1$, the set of 
matchings for vertex $1$ are $\MM(6,2)=\{\triangle_1,\triangle_2,\triangle_3,\triangle_4\}$ and 
$a_{6,2}=4$. 
For entry $a_{6,3}$ we have to consider matchings between triangles and the vertices 
$1,2$: $\MM(6,3)=\{(\triangle_1,\triangle_4), (\triangle_2,\triangle_4), (\triangle_3,\triangle_4)\}$ (first entry 
in the tuple is the triangle allocated to vertex 1, second entry the triangle allocated to vertex 2), $a_{6,3}=3$.  We compare this with $\MM(3,6)$ matching triangles to the vertices $4,5$ of 
\end{minipage}

\noindent
the polygon: 
$\MM(3,6)=\{(\triangle_3,\triangle_2),(\triangle_3,\triangle_1), (\triangle_2,\triangle_1)\}$ and $a_{3,6}=3$. 
\end{ex}

%
\subsection{A frieze determinant} 

As an immediate consequence of the geometric interpretation of Broline et al.\ of all frieze entries 
we recover the $n$-periodicity of friezes of order $n$ and we 
see that closed integral friezes are invariant under a glide reflection. 
We consider a fundamental 
domain for this glide reflection given by the entries in positions $(i,j)$ with 
$1\le i\le n-3$ (as in Figure~\ref{fig:order6-labels}). 
We include the $n-1$ entries equal to 1 and the $n$ entries equal to 0 above as well as 
the entry 1 below 
the two positions $(1,n-1)$ and $(2,n)$. 
These entries form a triangle and we take it as the upper triangular part of an $n\times n$ matrix, 
with a row of 0s on the diagonal. Reflecting along the 
diagonal, we create a symmetric matrix $M=M(T)$ whose entries are the matching numbers of a 
frieze pattern. This matrix depends on the chosen triangulation or on the corresponding 
matching numbers. But its determinant is independent of these choices, it only depends on the 
size of the polygon. 

\begin{theorem}\cite[Theorem 4]{BCI74}\label{thm:bci-2}
Let $M$ be the symmetric matrix corresponding to the frieze of a triangulation of an $n$-gon. 
Then $\det M= - (-2)^{n-2}$. 
\end{theorem}

The determinant result can be stated in terms of cluster variables. We associate the entries of a frieze 
of order $n$ with the diagonals of a convex $n$-gon. These in turn correspond to cluster variables 
$x_{ij}$, $1\le i<j-1\le n$, 
of a cluster algebra of type A$_{n-3}$, see \cite[\S 12.2]{FZ03}. 
The clusters of this cluster algebra are given by the triangulations of the $n$-gon, together 
with the frozen variables $\{x_{1,2},x_{2,3},\dots, x_{n-1,n},x_{n,1}\}$ of the edges of the polygon: 
$\{x_{ij}\mid (i,j) \mbox{ is a diagonal of $T$}\} \cup \{x_{i,i+1}\mid i=1,\dots, n\}$. 
The variables of the diagonals in the chosen triangulation are the {\bf initial cluster variables}. 
Every cluster variable of the associated cluster algebra is a Laurent polynomial in the initial variables, 
obtained through iterated mutations (a procedure similar to the diamond rule). Furthermore, 
the mutable cluster variables correspond bijectively to the diagonals in the $n$-gon. 
We put these variables in the positions as indicated by the endpoints of their diagonals. Then 
the cluster variables form a ``frieze'' $\FF=\FF(\underline{x})$ of elements of a ring of 
Laurent polynomials, in particular of an integral domain. 
This is a slight generalisation of the original definition of frieze patterns, replacing the two 
rows of 1s by the frozen variables $\{x_{1,2},x_{2,3},\dots, x_{n-1,n},x_{n,1}\}$ and naturally 
extending the labelling system with endpoints of diagonals in the $n$-gon. 
We note that $\FF(\underline{x})$ is tame. 

As before, in the resulting closed frieze there is a fundamental domain given by the $x_{ij}$ 
with $1\le i\le n-1$ and $i+1\le j\le n$. We use this fundamental domain as the upper triangular part 
of a square matrix $M=M(\underline{x})$ whose entries are cluster variables and 
make it symmetric around the diagonal of 0s. 

\begin{theorem}\cite[Theorem 3.2]{BM-determinant}\label{thm:bm-frieze}
Let $T$ be a triangulation of a convex $n$-gon with initial cluster 
$\underline{x}=\underline{x}(T)=\{x_{ij}\mid (ij)\ \in\ T\}\cup\{x_{i,i+1}\}_i$. 
Let $\FF(\underline{x})$ be the frieze of all the cluster variables obtained from this cluster and 
let $M(\underline{x})$ be the square matrix of $\FF(\underline{x})$. Then 
\[
{\rm det}\ M(\underline{x})= -(-2)^{n-2}x_{12}x_{23}\cdots x_{n-1,n}x_{n1}. 
\]
\end{theorem}

If we specialise the cluster variables $x_{ij}$ with $(ij)\in T$ and the frozen cluster variables to 1, 
the entries in the fundamental domain become the entries of the closed integral frieze pattern of $T$  
(\cite[Proposition 5.2]{CCh06}) and the frieze determinant from 
Theorem~\ref{thm:bm-frieze} specialises to the frieze determinant from Theorem~\ref{thm:bci-2}. 

\begin{rem} 
It would be interesting to see whether an analogous formula involving frozen variables 
can be found for cluster algebras in type 
D$_n$ which arise from triangulations of a punctured disk. 
\end{rem}

%
\section{Infinite friezes}\label{sec:infinite}
%

The theorem of Conway and Coxeter tells us how quiddity sequences of finite integral friezes arise. 
It is a natural question to ask which sequences yield infinite integral friezes patterns. 
An example which has already appeared in~\cite[Problem (16)]{Coco1} 
is the one on the right in Figure~\ref{fig:infinite}. In~\cite{T2015}, Manuela Tschabold showed that such patterns arise from 
triangulations of punctured disks. Let $S_n$ be a disk with $n$ marked points on the boundary and a point 
in the middle (the puncture). We use $n$ arcs to triangulate $S_n$: These arcs have as 
endpoints the vertices on the boundary or the point in the middle and they are pairwise non-crossing. 
The resulting regions have three sides in general; two of them may coincide. 
For each vertex $i$ of the boundary, we let $a_i$ be the number of 
connected components of the complement of the triangulation in a small neighbourhood of $i$. 

\begin{theorem}\cite[Theorem 3.6]{T2015}. 
Let $a_i$ be given by a triangulation of the punctured disk $S_n$. Then 
$(a_1,\dots, a_n)$ is the quiddity sequence of an infinite $n$-periodic integral frieze. 
\end{theorem}

Furthermore, as in the finite case every entry in the frieze of a triangulation 
of $S_n$ is a matching number between triangles and 
vertices (going to a suitable 
covering of the punctured disk), see~\cite[Section 4.5]{T2015}. 

\begin{ex}\label{ex:infinite-star}
Consider the star triangulation of a punctured disk with 4 vertices on the boundary on the left in Figure~\ref{fig:infinite}. 
At each vertex on the boundary, there are 2 triangles. We thus have $a_i=2$ for all $i$. 
The frieze associated to the quiddity sequence $(2,2,2,2)$ is on the right in Figure~\ref{fig:infinite}. 

\begin{figure}
\centering
\begin{minipage}{.4\textwidth}
\begin{tikzpicture}[scale=0.35]


    \fill[gray!40] (0,0) circle (4); 
    \draw[thick] (0,0) circle (115pt);
    \node[label=right:1,shape=circle,fill=black,scale=.4] (o1) at (0:4) {}; 
    \node[label=below:2,shape=circle,fill=black,scale=.4] (o2) at (270:4) {};
    \node[label=left:3,shape=circle,fill=black,scale=.4] (o3) at (180:4) {};
    \node[label=above:4,shape=circle,fill=black,scale=.4] (o4) at (90:4) {};
    \node[shape=circle,fill=black,scale=.4] (o5) at (0:0) {};

    
   \draw[thick, black] (o1) to (o5); 
   \draw[thick, black] (o2) to (o5);
   \draw[thick, black] (o3) to (o5);
   \draw[thick, black] (o4) to (o5);

\end{tikzpicture}
\end{minipage}
\begin{minipage}{.45\textwidth}
\[
\xymatrix@C=1em@R=.6em{
0 && 0 && 0 && 0 && 0   \\
& 1 && 1 && 1 && 1 && 1    \\
2 && 2 && 2 && 2 && 2 &   \\
& 3 && 3 && 3 && 3 && 3   \\
4 && 4 && 4 && 4 && 4 &   \\
 & \vdots && \vdots && \vdots && \vdots 
 }
\]
\end{minipage}
\caption{A triangulation of a punctured disk and its infinite integral frieze}\label{fig:infinite}
\end{figure}
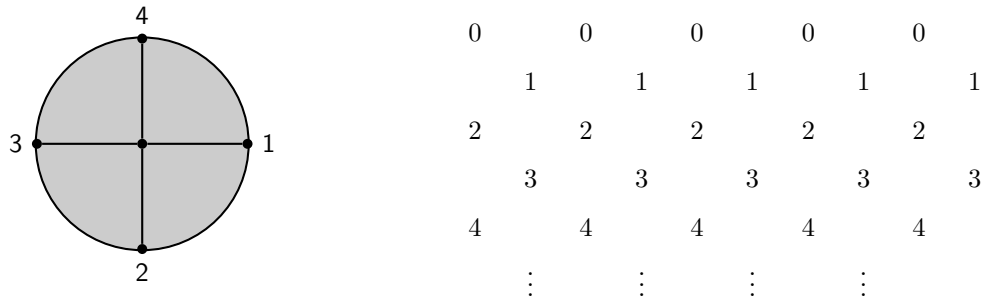
\end{ex}

The results of Tschabold show that there are infinite friezes for every period $n>0$ and that the infinite 
friezes constructed from triangulations of $S_n$ have a geometric interpretation as in the case of 
closed friezes. However, these are by far not the only types of infinite friezes. 
It can be shown that if $\mathcal F$ is an infinite $n$-periodic frieze with quiddity sequence 
$(a_1,\dots, a_n)$, then for any $b>0$, the $n$-tuple $(a_1+b,a_2,\dots, a_n)$ is also the 
quiddity sequence of an infinite $n$-periodic frieze, \cite[Theorem 2.2]{BPT}. 
Even if the original sequence can be realised by a triangulation of $S_n$, 
this is not true anymore for the new quiddity sequence. 

However, any quiddity sequence of an infinite frieze can be realised as the matching numbers of a 
triangulation of a surface. 

\begin{theorem}\cite[Theorem 4.6]{BPT}\label{thm:inf-annulus}
Each infinite periodic frieze comes from a triangulation of an annulus. All entries of the 
frieze are matching numbers between triangles and sets of vertices. 
\end{theorem}

\begin{rem}
To prove Theorem~\ref{thm:inf-annulus}, one can give an explicit construction of a triangulation 
for a quiddity sequence $(a_1,\dots, a_n)$ of an infinite frieze. This works as follows. 

(1) If there exists an entry 1 among the $a_i$, one reduces the quiddity sequence 
as in Remark~\ref{rem:cut-glue} until there are no 1s left. Since the frieze is infinite, after finitely 
many steps all entries in the quiddity sequence are $\ge 2$ (\cite[Remark 4.7]{BPT}). 

(2) Assume now that $a_i\ge 2$ for all $i$ and that there exist an entry $a_i>2$. Draw vertices $1,2,\dots, n$ 
on the outer boundary of an annulus. Then add  
$a_1-1$ starting segments of arcs at vertex $1$. Next, put $a_1-1$ vertices on the inner boundary, 
across  vertex 1 and connect the starting segments to these vertices on the inner boundary. 
In a second step, draw 
$a_2-1$ starting segments of arcs at vertex 2 . The first of them is connected to the last vertex 
drawn for the arcs at vertex 1. In addition, $a_2-2$ new vertices are drawn on the inner boundary. And so on. 
For vertex $n$, only $a_n-3$ new vertices are drawn on the inner boundary as the last of the $a_n-1$ 
arc at $n$ gets connected with the first vertex created for vertex 1. 

(3) If the quiddity sequence is $(2,2,\dots,2)$, we know that the frieze arises from a star triangulation of a punctured disk. 
We replace every arcs of the star triangulation by an arc starting at the outer boundary 
of an annulus and spiralling around a non-contractible 
closed curve in this annulus (all in the same direction). See Example~\ref{ex:spiraling}. 
\end{rem}

\begin{ex}\label{ex:spiraling} 
The triangulation of an annulus giving rise to the trivial quiddity sequence $(2,2,2,2)$ is given 
by arcs spiraling around a non-contractible curve in the annulus as shown on the left hand side 
in Figure~\ref{fig:all-twos}. The triangulation yielding the quiddity sequence $(3,4,2,4)$ is 
shown in the right.

\begin{figure}[ht]
\centering
\begin{minipage}{.4\textwidth}
\begin{tikzpicture}[scale=.4]
\draw[thick,fill=gray!40, radius=4cm] (0,0) circle; 
\draw[thick,fill=white!10, radius=1cm] (0,0) circle; 
\draw[thick, radius=1.9cm, dashed] (0,0) circle;

\node[inner sep=1.5pt, fill, circle] (3) at (0:4) {};
\node[inner sep=1.5pt, fill, circle] (4) at (90:4) {};
\node[inner sep=1.5pt, fill, circle] (1) at (180:4) {};
\node[inner sep=1.5pt, fill, circle] (2) at (270:4) {};
\node[left] (1txt) at (1) {$3$};
\node[below] (2txt) at (2) {$2$};
\node[right] (3txt) at (3) {$1$};
\node[above] (4txt) at (4) {$4$};

\begin{scope}[thick, rounded corners=10pt]
\draw (1) -- ++(0:2) .. controls (210:2.1) and (240:2.2) .. (270:1.92); 
\draw (2) -- ++(90:2) .. controls (300:2.1) and (330:2.2) .. (360:1.92); 
\draw (3) -- ++(180:2) .. controls (30:2.1) and (60:2.2) .. (90:1.92); 
\draw (4) -- ++(270:2) .. controls (120:2.1) and (150:2.2) .. (180:1.92); 
\end{scope}
\end{tikzpicture}
%
%
\end{minipage}
\begin{minipage}{.4\textwidth}
\begin{tikzpicture}[scale=.4]
\draw[thick,fill=gray!40, radius=4cm] (0,0) circle; 
\draw[thick,fill=white!10, radius=2cm] (0,0) circle; 
\node[inner sep=1.5pt, fill, circle] (1) at (0:4) {};

\node[inner sep=1.5pt, fill, circle] (2) at (90:4) {};

\node[inner sep=1.5pt, fill, circle] (3) at (180:4) {};

\node[inner sep=1.5pt, fill, circle] (4) at (270:4) {};

\node[right] (1txt) at (1) {$1$};
\node[above] (2txt) at (2) {$4$};
\node[left] (3txt) at (3) {$3$};
\node[below] (4txt) at (4) {$2$};
\node[inner sep=1.5pt, fill, circle] (1-bar) at (0:2) {};
\node[inner sep=1.5pt, fill, circle] (2-bar) at (70:2) {};
\node[inner sep=1.5pt, fill, circle] (3-bar) at (130:2) {};
\node[inner sep=1.5pt, fill, circle] (4-bar) at (190:2) {};
\node[inner sep=1.5pt, fill, circle] (5-bar) at (290:2) {};

\begin{scope}[thick]
\draw	 (3) to (4-bar);
\draw	 (4) to[bend left] (4-bar);
\draw	 (4) to (5-bar);
\draw	 (4) to[bend right] (1-bar);
\draw	 (1) to[bend right] (2-bar);
\draw	 (2) to (2-bar);
\draw	 (2) to (3-bar);
\draw	 (2) .. controls (-3,2) and (-2.5,1) .. (4-bar);
\end{scope}
\draw[thick] (1) to (1-bar);

\begin{scope}[even odd rule, thick]
\clip (1) circle (.8cm) (2) circle (.8cm) (3) circle (.8cm) (4) circle (.8cm);
\draw	 (3) to (4-bar);
\draw	 (4) to[bend left] (4-bar);
\draw	 (4) to (5-bar);
\draw	 (4) to[bend right] (1-bar);
\draw	 (1) to[bend right] (2-bar);
\draw	 (2) to (2-bar);
\draw	 (2) to (3-bar);
\draw	 (2) .. controls (-3,2) and (-2.5,1) .. (4-bar);
\draw[thick] (1) to (1-bar);
\end{scope}

\end{tikzpicture}
\end{minipage}
\caption{Triangulations for the quiddity sequences $(2,2,2,2)$ and $(3,4,2,4)$.} 
\label{fig:all-twos}
\end{figure}
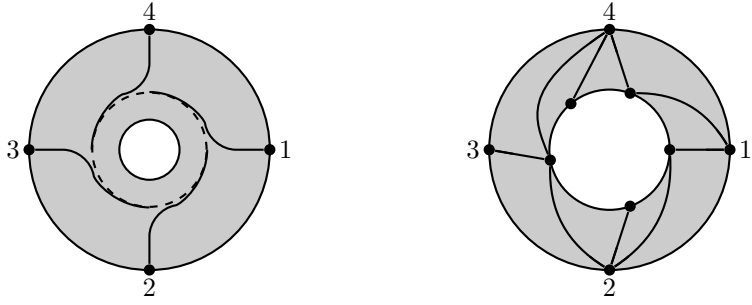
\end{ex}

%
\subsection{Growth of frieze patterns} 

Clearly, in any closed frieze, the entries along a diagonal go from 0 over 1 to positive integers and then back to 1 and 0. 
If an infinite frieze arises from a triangulation of a punctured disk, it can be shown that 
the entries in any diagonal form 
$n$ arithmetic progressions, \cite[Proposition 3.12]{T2015}. 
An illustration is the frieze pattern in Example~\ref{ex:infinite-star}. 
As Example~\ref{ex:fast} illustrates, the entries in an infinite frieze can grow much faster. 
In fact, we will see that in general, the numbers grow exponentially. 

\begin{ex}\label{ex:fast}
The entries in frieze patterns arising from triangulations of annuli grow very fast. An illustration for this is the 
frieze with quiddity sequence $(3,4,2,4)$ as we see here. The entries in the 5th non-trivial row are already 3digit numbers. 
\[
\xymatrix@C=.7em@R=.6em{
0 && 0 && 0 && 0 && 0 && 0  \\
& 1\ar@{..}[dddddd]   && 1 && 1 && 1 && 1\ar@{..}[dddddd]   \\
4 && 3 && 4 && 2 && 4 && 3   \\
& 11 && 11 && 7 && 7 && 11   \\
19 && 40 && 19 && 24 && 19 && 40  \\
& 69 && 69 && 65 && 65 && 69   \\
236 && 119 && 236 && 176 && 236 && 119  \\
 &  && \vdots && \vdots && \vdots && 
}
\]
\end{ex}

Let $\FF$ be an integral frieze - finite or infinite. Let $n_0$ be the minimal period it has. If $\FF$ is finite of order $n$, it arises from a triangulation 
of an $n$-gon. Then $n_0$ is equal to $n/2$ 
if the triangulation is invariant under rotation by 180 degrees, 
equal to $n/3$ if the triangulation is invariant under a rotation by 60 degrees. It is equal to $n$ otherwise. 

\begin{ex}\label{ex:s-finite}
We consider the closed friezes of small order. If $\FF$ is of order 4, its quiddity sequence is $(1,2,1,2)$ or $(2,1,2,1)$. In both cases, 
the smallest period is 2. If $\FF$ has order 5, it arises from a triangulation of a pentagon and 
since no such triangulation has a non-trivial rotational symmetry, the smallest period is 5. 
If $\FF$ has order 6, it comes from a triangulation of a hexagon. Among them, there are the triangulations with rotational symmetry 
by 180 degrees, e.g. with diagonals $(1,3)$, $(3,6)$, $(4,6)$ and the triangulations with 3-fold symmetry, e.g. by the diagonals 
$(1,3)$, $(1,5)$ and $(3,5)$. In the former case, the frieze has period 3, in the latter case, it is 2-periodic. 

Here we show the 5-periodic frieze with quiddity sequence $(1,2,2,1,3)$ and the 2-periodic frieze with quiddity sequence $(1,3,1,3,1,3)$. 
For later use, we add rows of $-1$s at both ends.
\[
\xymatrix@C=.1em@R=.3em{
 & -1 && -1 && -1 && -1 && -1 &&\ \ &&  -1 && -1 && -1 && -1 && -1&& -1\\
0 && 0 && 0 && 0 && 0 && 0  &&&& 0 && 0 && 0 && 0 && 0 && 0\\
& 1  && 1 && 1 && 1 && 1  &&&&  1  && 1 && 1 && 1 && 1 && 1\\
1 && 2 && 2 && 1 && 3 && 1 &&&&  1 && 3 && 1 && 3 && 1  && 3 \\
& 1 && 3 && 1 && 2 && 2   &&&& 2 && 2 && 2 && 2 && 2 && 2 \\
1 && 1 && 1 && 1 && 1 && 1  &&&& 3 && 1 && 3 && 1 && 3 && 1 \\
& 0 && 0 && 0 && 0 && 0  &&&& 1  && 1 && 1 && 1 && 1 && 1 \\
-1 && -1 && -1 && -1 && -1 && -1 &&&&  0 && 0 && 0 && 0 && 0 && 0 \\ 
& && && && && &&  &&      -1 && -1 && -1 && -1 && -1 && -1 
}
\]
\end{ex}

We recall that the elements 
in a frieze are indexed so that the entries in the quiddity sequence are of the form $a_{i,i+2}$ and going down a diagonal 
NW-SE, the first coordinate is fixed while the second increases; going one position to the right, both coordinates increase by 1 
(Definition~\ref{def:frieze}). We extend this and denote the entries in the first rows of 1s by $a_{i,i+1}$, the entries in the rows of 0s above by 
$a_{i,i}$ and the entries $-1$ above them by $a_{i,i-1}$. 
When drawing a frieze in the plane, we say that the (first) row containing the quiddity sequence is row 1, the rows above are called 
rows 0, -1 and -2 (the latter is the row of $-1$s written above the frieze). If the frieze is closed of 
order $n$, 
the frieze ends with a row of 1s (row $n-2$) followed by a row of 0s (row $n-1$) and we extend 
this by the additional row of $-1$s (row $n$). 

We have everything ready to state a remarkable property of integral friezes: 
\begin{theorem}\cite[Theorem 2.2]{BFPT} 
Let $\FF$ be an integral frieze. Assume that $\FF$ is $m$-periodic. Then for any $k\ge 0$, 
the differences of the entries in rows $km$ and $km-2$ is constant: 
$a_{i,i+km+1}-a_{i+1,i+km}=a_{j,j+km+1}-a_{j+1,j+km}$ for all $i,j$.  
\end{theorem}

There are two extreme cases where the statement is clearly true: If $k=0$, it concerns 
the difference of the entries in rows $0$ and 
$-2$, this is always 2. Secondly, 
if $\FF$ is closed of order $n$, it is in particular $n$-periodic. For closed friezes of order $n$, we reduce the coordinates of the entries 
modulo $n$ as indicated in Figure~\ref{fig:order6-labels}. So if $k=1$ and $m=n$ 
we have $a_{i,i+n+1}=a_{i,i+1}=-1$, 
$a_{i+1,i+n}=a_{i+1,i}=1$ for all $i$, with constant difference $-2$. 

These constant differences can be expressed using Chebychev polynomials of the first kind. More precisely: 
For infinite friezes, the constant differences of entries occur every $n_0$ rows, if $n_0$ is the smallest period. We define $s:= a_{i,i+n_0+1}-a_{i+1,i+n_0}$ to be the 
first of these differences (apart from the one for $k=0$). 
In the frieze pattern of Example~\ref{ex:infinite-star}, $s=2$. In 
Example~\ref{ex:fast}, we have $s=58$. 

Note that if the 
frieze arises from a triangulation of a punctured disk, all these constants are equal to 2 (the first case to check here is the frieze we have 
seen before, from Figure~\ref{fig:infinite}). 
In all other infinite friezes, these differences grow fast. 
The sequence 
$(s_k)_{k\ge 0}$ of these constant differences satisfies 
the recurrence 
\[
s_{k+2}=s_1s_{k+1}-s_k
\]
Explicitly, in terms of $s=s_1$, we have 
\[
s_k=s^k+k\sum_{l=1}^{\lfloor k/2\rfloor}(-1)^l\frac{1}{k-l}{k-l \choose l} s^{k-2l}
\]
for $k\ge 1$. 
Furthermore, the sequence $(s_k)_{k\ge 0}$ grows asymptotically 
exponentially if and only if $s>1$. In this case, its 
growth rate is $\frac{1+\sqrt{s^2-4}}{2}$, see \cite[\S 4]{BFPT} for details. 
As $s$ governs the rate at which the entries in $\FF$ grow, it is called the {\bf growth coefficient} 
of the frieze. 

\begin{rem}
(1) 
Any triangulated annulus with marked points on both boundaries 
gives rise to two frieze patterns, as both boundary components determine a quiddity 
sequence. These sequences are different in general. However, one can show that their 
growth coefficients are the same, \cite[Theorem 3.4]{BFPT}. 
In the figure on the right in Example~\ref{ex:fast}, the quiddity 
sequence of the inner boundary is $(3,2,4,2,3)$, giving a 5-periodic frieze. 

(2) Triangulations of annuli give rise to cluster tilting objects in a cluster category of type 
$\widetilde{A}$. In this algebraic interpretation, the common growth behaviour of the two associated 
friezes can be viewed as common behaviour of tubes in the associated category and 
interpreted as generalised numbers of submodules of indecomposable 
objects, \cite[Section 4.4]{BCJKT20}. However there are examples cluster categories with collections 
of tubes where the growth coefficients differ across the different tubes. 
\end{rem}

\bibliographystyle{plain}
\bibliography{biblio}

\end{document}